\input ./style/arxiv-general.cfg
\documentclass[MSNbibl,number,citesort,secthm,seceqn,dvips]{arxbj}
\makeatletter
   \@ifpackageloaded{graphicx}{}{\usepackage{graphicx}}
\makeatother
\usepackage{upgreek}

\volume{22}
\issue{3}
\pubyear{2016}
\firstpage{1598}
\lastpage{1616}
\doi{10.3150/15-BEJ705}
\docsubty{FLA}

\makeatletter
\newtheorem{theorem}{Theorem}[section]
\newtheorem{proposition}[theorem]{Proposition}
\newremark{remark}[theorem]{Remark}
\newremark{example}[theorem]{Example}
\makeatother

\begin{document}
\begin{frontmatter}

\title{$L^p$-Wasserstein distance for stochastic differential equations
driven by L\'{e}vy~processes}
\runtitle{$L^p$-Wasserstein distance for stochastic differential
equations driven by L\'{e}vy processes}

\begin{aug}
\author{\inits{J.}\fnms{Jian}~\snm{Wang}\corref{}\ead
[label=e1]{jianwang@fjnu.edu.cn}}
\address{School of Mathematics and Computer Science, Fujian Normal
University,
Fuzhou 350007, China. \printead{e1}}
\end{aug}

%
\received{\smonth{8} \syear{2014}}
%
\revised{\smonth{1} \syear{2015}}

%
\begin{abstract}
Coupling by reflection mixed with synchronous coupling is constructed
for a class of stochastic differential equations (SDEs) driven by L\'{e}vy noises.
As an application, we establish the exponential contractivity of the
associated semigroups $(P_t)_{t\ge0}$ with respect to the standard
$L^p$-Wasserstein distance for all $p\in[1,\infty)$. In particular,
consider the following SDE:
\[
\mathrm{d}X_t=\mathrm{d}Z_t+b(X_t)\,
\mathrm{d}t,
\]
where $(Z_t)_{t\ge0}$ is a symmetric $\alpha$-stable
process on $\mathbb{R}^d$ with $\alpha\in(1,2)$. We show that if the
drift term
$b$ satisfies that for any $x,y\in\mathbb{R}^d$,
\[
\bigl\langle b(x)-b(y),x-y\bigr\rangle\le\cases{ K_1|x-y|^2,\hspace*{8pt}
\qquad |x-y|\le L_0;\vspace*{3pt}
\cr
-K_2|x-y|^\theta,
\qquad |x-y|> L_0}
\]
holds with some positive constants $K_1$, $K_2$, $L_0>0$ and
$\theta\ge2$, then there is a constant $\lambda:=\lambda(\theta,
K_1,K_2, L_0)>0$ such that for
all $p\in[1,\infty)$, $t>0$ and $x,y\in\mathbb{R}^d$,
\[
W_p(\delta_x P_t, \delta_y
P_t)\le C(p,\theta, K_1,K_2,L_0)\mathrm{e}^{-\lambda t/p}
\biggl[\frac{|x-y|^{1/p}\vee|x-y|}{1+|x-y|{\mathbf1}_{(1,\infty
)\times
(2,\infty)}(t,\theta)} \biggr].
\]
\end{abstract}

%
\begin{keyword}
\kwd{coupling by reflection}
\kwd{exponential contractivity}
\kwd{$L^p$-Wasserstein distance}
\kwd{stochastic differential equation driven by L\'{e}vy noise}
\kwd{symmetric stable process}
\end{keyword}
\end{frontmatter}

\section{Introduction}\label{section1}

In this paper, we consider the following stochastic differential
equation (SDE) driven by L\'{e}vy noises:
%
\begin{equation}
\label{Ito-SDE} \mathrm{d}X_t=\mathrm{d}Z_t+b(X_t)
\,\mathrm{d}t,
\end{equation}
where $(Z_t)_{t\ge0}$ is a $d$-dimensional L\'{e}vy process, and
$b \dvtx\mathbb{R}^d\to\mathbb{R}^d$ is a continuous vector field
such that for any
$x,y\in\mathbb{R}^d$,
\[
\bigl\langle b(x)-b(y),x-y\bigr\rangle\le C|x-y|^2
\]
holds for some constant $C>0$. It is a standard result that in this
case the SDE
(\ref{Ito-SDE}) enjoys the unique strong solution.

Denote by $(P_t)_{t\ge0}$ the semigroup associated to
(\ref{Ito-SDE}). If the initial value $X_0$ is distributed as $\mu$,
then for any $t>0$, the distribution of $X_t$ is $\mu P_t$. We are
concerned with the exponential contractivity of
the map $\mu\mapsto\mu P_t$ with respect to the standard $L^p$-Wasserstein
distance $W_p$ for all $p\geq1$.
Given two probability measures $\mu$ and $\nu$ on $\mathbb{R}^d$,
the standard
$L^p$-Wasserstein distance
$W_p$ for all $p\in[1,\infty)$ (with respect to the Euclidean norm
$|\cdot|$) is given by
\[
W_p(\mu,\nu)=\inf_{\Pi\in\mathcal{C}(\mu,\nu)} \biggl(\int
_{\mathbb{R}^d\times\mathbb{R}^d} |x-y|^p\,\mathrm{d}\Pi(x,y)
\biggr)^{1/p}.
\]
Equipped with $W_p$, the totality $\mathcal P_p(\mathbb{R}^d)$ of probability
measures having finite moment of order $p$
becomes a complete metric space.

The following result is well known.

\begin{theorem}\label{theorem1}
Suppose that there exists a constant $K>0$ such that
%
\begin{equation}
\label{diss-1}
\bigl\langle b(x)-b(y),x-y\bigr\rangle\le-K|x-y|^2\qquad
\mbox{for all } x,y\in \mathbb{R}^d.
\end{equation}
Then, for any $p\ge1$ and $t>0$,
%
\begin{equation}
\label{w2}
W_p(\mu P_t,\nu P_t)\le
\mathrm{e}^{-Kt} W_p(\mu,\nu)\qquad\mbox{for all } \mu ,\nu\in\mathcal
P_p\bigl(\mathbb{R}^d\bigr).
\end{equation}
\end{theorem}

The proof of this result is quite straightforward, by simply using
the synchronous coupling, which is also called the basic coupling or
the coupling of marching
soldiers in the literature (see, e.g., \cite{CL}, Definition~2.4 and
\cite{Chen}, Example~2.16). The
reader can refer to \cite{BGG}, page~2432, the proof of Theorem~\ref
{theorem1} for the case of diffusion processes. (\ref{diss-1}) is
the so-called uniformly dissipative condition, which seems to be a
limit in applications.
For diffusion processes, it follows from \cite{RS}, Theorem~1, or \cite{BGG}, Remark~3.6 (also see \cite{CG}, Theorem~3.6) that (\ref{w2}) holds
for any probability measures $\mu$ and $\nu$ if and only if (\ref
{diss-1}) holds for all $x$, $y\in\mathbb{R}^d$.
The first breakthrough to get rid of such restrictive condition in this
direction for $L^1$-Wasserstein distance $W_1$ was done recently by
Eberle in \cite{Eberle,Eberle11},
at the price of multiplying a constant $C\geq1$ on the right-hand side
of (\ref{w2}). See \cite{Eberle}, Corollary~2.3,  for more details, and
\cite{LW}, Theorem~1.3,  for related developments on $L^p$-Wasserstein
distance $W_p$ with all $p\in[1,\infty)$ on this topic. However, the
corresponding result for SDEs driven by L\'{e}vy noises is not
available yet now. Indeed, we will see later that in this case we need
a completely different idea for the construction of the coupling
processes, and a new approach by using the coupling argument, in
particular the more delicate choice of auxiliary functions.

Throughout this paper, we assume that the driving L\'{e}vy process
has a symmetric $\alpha$-stable process as a component. That is, let
$\nu$ be the L\'{e}vy measure of the process $(Z_t)_{t\ge0}$, then
\[
\nu(\mathrm{d}z)\ge\frac{C_{d,\alpha}}{|z|^{d+\alpha}}\,\mathrm{d}z,
\]
where
$C_{d,\alpha}=2^\alpha\Gamma((d+\alpha)/2)\uppi^{-d/2}
|\Gamma(-\alpha/2)|^{-1}$ is a constant associated with the L\'{e}vy
measure of a symmetric $\alpha$-stable process or fractional
Laplacian, that is,
\[
-(-\Delta)^{\alpha/2} f(x)=\int \bigl(f(x+z)-f(x)-\bigl\langle\nabla f(x) ,
z\bigr\rangle\mathbf{1}_{\{|z|\le1\}} \bigr) \frac{C_{d,\alpha
}}{|z|^{d+\alpha}}\,\mathrm{d}z.
\]
Denote by $\omega_d=2\uppi^{d/2}/\Gamma(d/2)$ the surface measure of
the unit sphere in $\mathbb{R}^d$. Our main contribution of this paper is
as follows.

\begin{theorem}\label{main-result}
Assume that for any $x,y\in\mathbb{R}^d$,
%
\begin{equation}
\label{main-result-1}
\bigl\langle b(x)-b(y),x-y\bigr\rangle\le %
\cases{
K_1|x-y|^2,&\quad$|x-y|\le L_0$;
\vspace*{3pt}
\cr
-K_2|x-y|^\theta,&\quad$|x-y|>
L_0$}
\end{equation}
holds with some positive constants $K_1$, $K_2$, $L_0>0$ and
$\theta\ge2$.
Then for all $\alpha\in(1,2)$ or for all $\alpha\in(0,1]$ with
%
\begin{equation}
\label{main-result-1-1}
\frac{\alpha C_{d,\alpha
}\omega
_d3^{\alpha-1}}{8(2-\alpha)d}> K_1 L_0^\alpha,
\end{equation}
there exists a positive constant
$\lambda:=\lambda(\theta,K_1, K_2, L_0)>0$, such that for any $p\ge
1$ the
following two statements hold:
\begin{longlist}[(ii)]
\item[(i)] if $\theta=2$, then for all
$x$, $y\in\mathbb{R}^d$ and any $t>0$,
%
\begin{equation}
\label{main-result-2}
W_{p}(\delta_x P_t,
\delta_y P_t)\leq C\mathrm{e}^{-\lambda t/p} \bigl(|x-y|^{1/p}
\vee|x-y| \bigr);
\end{equation}

\item[(ii)] if $\theta>2$, then for all $x$, $y\in\mathbb{R}^d$
and any $t>0$,
%
\begin{equation}
\label{main-result-3} %
W_{p}(\delta_x
P_t,\delta_y P_t)\leq C\mathrm{e}^{-\lambda
t/p}
\biggl[\frac{|x-y|^{1/p}\vee|x-y|}{1+|x-y|{\mathbf1}_{(1,\infty
)}(t)} \biggr],
\end{equation}
\end{longlist}
where $C>0$ is a positive constant depending on $\theta$, $K_1$, $K_2$,
$L_0$ and $p$.
\end{theorem}

Theorem~\ref{main-result} above does provide new conditions on the
drift term $b$ such that the associated semigroup $(P_t)_{t\ge0}$ is
exponentially contractive with respect to the $L^p$-Wasserstein
distance $W_p$ for all $p\ge1$. In particular, when
$\alpha\in(1,2)$, the conclusion of Theorem~\ref{main-result} is the
same as that of~\cite{LW}, Theorem~1.3, for diffusion processes;
while for $\alpha\in(0,1]$ we need the restrictive condition
(\ref{main-result-1-1}); see Remark~\ref{main-1-1} for a further
comment. Indeed, (\ref{main-result-1-1}) is natural in the sense
that, when $\alpha\in(0,1]$ the drift term plays the dominant
role or the same role (just in case that $\alpha=1$) for the
behavior of SDEs driven by symmetric $\alpha$-stable processes,
see, for example, \cite{BJ,Chen-Song} for (Dirichlet) heat kernel estimates
and \cite{Wang-Wang} for dimensional free Harnack inequalities on
this topic. Similarly, in considering the exponential contractivity
of SDE (\ref{Ito-SDE}), we need (\ref{main-result-1-1}) to control
the locally non-dissipative part of the drift term. Note that
(\ref{main-result-1-1}) holds true when $K_1, L_0$ are small enough.

To show the power of Theorem~\ref{main-result}, we consider the
following example about the SDE driven by symmetric $\alpha$-stable
processes with $\alpha\in(0,2)$, which yields the exponential
contractivity of the
semigroup $(P_t)_{t\ge0}$ with respect to the $L^p$-Wasserstein
distance $W_p$ $(p\ge1)$ for super-convex potentials.

\begin{example}\label{1-example-2}
Let $(Z_t)_{t\ge0}$ be a symmetric $\alpha$-stable process in
$\mathbb{R}^d$
with $\alpha\in(0,2)$, and $b(x)=\nabla V(x)$ with
$V(x)=-|x|^{2\beta}$ and $\beta>1$. Then there exists a constant
$\lambda:=\lambda(\alpha,\beta)>0$ such that for all $p\ge1$, $x$,
$y\in\mathbb{R}^d$ and $t>0$,
\[
W_{p}(\delta_x P_t,\delta_y
P_t)\leq C(\alpha,\beta,p) \mathrm{e}^{-\lambda t/p
} \biggl[\frac{|x-y|^{1/p}\vee|x-y|}{1+|x-y|{\mathbf1}_{(1,\infty
)}(t)}
\biggr].
\]
\end{example}

Note that the uniformly dissipative condition (\ref{diss-1}) fails for
Example~\ref{1-example-2}; see, for example, (\ref{cg-1}) below.
That is, one cannot deduce directly from Theorem~\ref{theorem1} the exponential
contractivity with respect to the $L^{p}$-Wasserstein distance $W_{p}$
for all $p\in[1,\infty)$.

The remainder of this paper is arranged as follows. In the next
section, we will present the coupling by reflection mixed with the
synchronous coupling for the SDE (\ref{Ito-SDE}) driven by L\'{e}vy
noise, and also prove the existence of coupling process associated with
this coupling (operator). Section~\ref{section3} is mainly devoted to
the proof of Theorem~\ref{main-result}. For this, we need more delicate
choice of auxiliary functions and some key estimates for them, which
are different between $\alpha\in(1,2)$ and $\alpha\in(0,1]$. The sketch
of the proof of Example~\ref{1-example-2} is also given here.

\section{Coupling operator and coupling process for SDEs with
jumps}

\subsection{Coupling by reflection and synchronous coupling}
It is easy
to see that the generator of the process $(X_t)_{t\ge0}$ acting on
$C_b^2(\mathbb{R}^d)$ is
%
\begin{equation}
\label{ope-0}
Lf(x)=\int \bigl(f(x+z)-f(x)-\bigl\langle \nabla f(x) , z\bigr
\rangle\mathbf{1}_{\{|z|\le1\}} \bigr) \nu(\mathrm{d}z)+ \bigl\langle b(x),
\nabla f(x)\bigr\rangle.
\end{equation}
In this part, we construct a coupling operator for the
generator $L$ above. For any $x$, $y$ and
$z\in\mathbb{R}^d$, we set
\[
\varphi_{x,y}(z):= %
\cases{ \displaystyle z-\frac{2\langle x-y, z\rangle}{|x-y|^2}(x-y), &
\quad $x\neq y$;
\vspace*{3pt}\cr
-z, & \quad $x=y$.}
\]
It is clear that $\varphi_{x,y}\dvtx \mathbb{R}^d\to\mathbb{R}^d$ has
the following three
properties:
\begin{longlist}[(A3)]
\item[(A1)] $\varphi_{x,y}(z)=\varphi_{y,x}(z)$ and $\varphi
^2_{x,y}(z)=z$, that is,
$\varphi_{x,y}^{-1}(z)=\varphi_{x,y}(z)$;
\item[(A2)] $|\varphi_{x,y}(z)|=|z|$;
\item[(A3)] $(z-\varphi_{x,y}(z))\, /\!\!/ \, (x-y)$ and $(z+\varphi
_{x,y}(z))\perp(x-y)$.
\end{longlist}
Next, for any $f\in C_b^2(\mathbb{R}^{2d})$, let
\[
\nabla_xf(x,y):= \biggl(\frac{\partial f(x,y)}{\partial x_i} \biggr)_{1\le
i\le d},\qquad
\nabla_yf(x,y):= \biggl(\frac{\partial f(x,y)}{\partial y_i} \biggr)_{1\le i\le d}.
\]

Now, let $L_0$ be the
constant appearing in (\ref{main-result-1}). We will split the
construction of the coupling operator into two parts, according to
$x,y\in\mathbb{R}^d$ with $|x-y|\le L_0$ or with $|x-y|> L_0$. First,
for any
$f\in C_b^2(\mathbb{R}^{2d})$ and $x,y\in\mathbb{R}^d$ with
$|x-y|\le L_0$, we define
\begin{eqnarray*}
\widetilde{L}f(x,y)&:=& \frac{1}{2} \biggl[\int_{ \{|z|\le
a{|x-y|} \}}
\bigl(f\bigl(x+z,y+\varphi_{x,y}(z)\bigr)-f(x,y)-\bigl\langle \nabla
_xf(x,y), z\bigr\rangle \mathbf{1}_{\{|z|\le1\}}
\\
&&\hspace*{44pt}\qquad{}-\bigl\langle\nabla_yf(x,y),\varphi_{x,y}(z)
\bigr\rangle \mathbf {1}_{\{|z|\le1\}} \bigr) \frac{C_{d,\alpha}}{|z|^{d+\alpha}}\,\mathrm{d}z
\\
&&\hspace*{12pt}{}+\int_{ \{|z|\le a{|x-y|} \}} \bigl( f\bigl(x+\varphi_{x,y}(z),y+z
\bigr)-f(x,y) -\bigl\langle\nabla_yf(x,y), z\bigr\rangle
\mathbf{1}_{\{|z|\le1\}}
\\
&&\hspace*{75pt}{}-\bigl\langle\nabla_xf(x,y), \varphi_{x,y}(z)
\bigr\rangle \mathbf{1}_{\{|z|\le1\}} \bigr) \frac{C_{d,\alpha}}{|z|^{d+\alpha}}\,\mathrm{d}z
\biggr]
\\
&&{}+\int_{ \{|z|\le a{|x-y|} \}} \bigl(f(x+z,y+z)-f(x,y)-\bigl\langle
\nabla_xf(x,y), z\bigr\rangle \mathbf{1}_{\{
|z|\le1\}}
\\
&&\hspace*{63pt}{}-\bigl\langle\nabla_yf(x,y), z\bigr\rangle
\mathbf{1}_{\{|z|\le
1\}} \bigr) \biggl(\nu(\mathrm{d} z)-\frac{C_{d,\alpha}}{|z|^{d+\alpha}}\,
\mathrm{d}z \biggr)
\\
&&{}+\int_{ \{|z|>a{|x-y|} \}} \bigl(f(x+z,y+z)-f(x,y)-\bigl\langle
\nabla_xf(x,y), z\bigr\rangle \mathbf{1}_{\{|z|\le1\}}
\\
&&\hspace*{63pt}{}-\bigl\langle\nabla_yf(x,y), z\bigr\rangle
\mathbf{1}_{\{|z|\le1\}} \bigr) \nu(\mathrm{d}z)
\\
&&{}+\bigl\langle\nabla_x f(x,y), b(x)\bigr\rangle+ \bigl
\langle\nabla_y f(x,y), b(y)\bigr\rangle,
\end{eqnarray*}
where $a\in(0,1/2)$ is a constant determined by
later.

On the other hand, for any $f\in C_b^2(\mathbb{R}^{2d})$ and $x,y\in
\mathbb{R}^d$
with $|x-y|> L_0$, we define
\begin{eqnarray*}
\widetilde{L}f(x,y)&:=&\int \bigl(f(x+z,y+z)-f(x,y)-\bigl\langle\nabla
_xf(x,y), z\bigr\rangle \mathbf{1}_{\{|z|\le1\}}
\\
&&\hspace*{13pt}{}-\bigl\langle\nabla_yf(x,y), z\bigr\rangle
\mathbf{1}_{\{|z|\le1\}} \bigr) \nu(\mathrm{d}z)
\\
&&{}+\bigl\langle\nabla_x f(x,y), b(x)\bigr\rangle+ \bigl
\langle\nabla_y f(x,y), b(y)\bigr\rangle.
\end{eqnarray*}
We can conclude the following.

\begin{proposition} The operator $ \widetilde{L}$ defined by above is the
coupling operator of the operator $L$ given by
(\ref{ope-0}).
\end{proposition}

\begin{pf} Since $\widetilde{L}$ is a linear operator, it suffices
to verify that
%
\begin{equation}
\label{coup-proof}\widetilde{L} f(x)=Lf(x),\qquad f\in C_b^2\bigl(
\mathbb{R}^d\bigr),
\end{equation}
where, on the left-hand side, $f$ is regarded as a bivariate
function on $\mathbb{R}^{2d}$, that is, $f(x)=f(x,y)$ for all $x, y\in
\mathbb{R}^d$.

For any $x,y\in\mathbb{R}^d$ with $|x-y|> L_0$, it is trivial to see
that (\ref{coup-proof}) holds true, and so we only need to verify that
for $x,y\in
\mathbb{R}^d$ with $|x-y|\le L_0$.
First, we have
\begin{eqnarray*}
\widetilde{L}f(x) &=&\frac{1}{2} \biggl[\int_{ \{|z|\le a{|x-y|} \}}
\bigl(f(x+z)-f(x)-\bigl\langle\nabla f(x), z\bigr\rangle \mathbf{1}_{\{|z|\le
1\}}
\bigr)\frac
{C_{d,\alpha}}{|z|^{d+\alpha}}\,\mathrm{d}z
\\
&&\hspace*{12pt}{}+\int_{ \{|z|\le a{|x-y|} \}} \bigl( f\bigl(x+\varphi
_{x,y}(z)\bigr)-f(x)-\bigl\langle\nabla f(x), \varphi_{x,y}(z)
\bigr\rangle \mathbf {1}_{\{|z|\le1\}} \bigr)
\frac{C_{d,\alpha}}{|z|^{d+\alpha}}\,\mathrm {d}z \biggr]
\\
&&{}+\int_{ \{|z|\le a{|x-y|} \}} \bigl(f(x+z)-f(x)-\bigl\langle \nabla
f(x), z\bigr\rangle \mathbf{1}_{\{|z|\le1\}} \bigr) \biggl(\nu (\mathrm{d}z)-
\frac
{C_{d,\alpha}}{|z|^{d+\alpha}} \,\mathrm{d}z \biggr)
\\
&&{}+\int_{ \{|z|> a{|x-y|} \}} \bigl(f(x+z)-f(x)-\bigl\langle \nabla
f(x), z\bigr\rangle \mathbf{1}_{\{|z|\le1\}} \bigr) \nu(\mathrm {d}z)
\\
&&{}+\bigl\langle b(x),\nabla f(x)\bigr\rangle.
\end{eqnarray*}
By (A1) and (A2), we know that the measure
$\frac{C_{d,\alpha}}{|z|^{d+\alpha}} \,\mathrm{d}z$ is
invariant under the transformation $z\mapsto\varphi_{x,y}(z)$.
This, along with (A2) and the equality above, leads to
\begin{eqnarray*}
\widetilde{L} f(x) &=&\int_{ \{|z|\le a{|x-y|} \}} \bigl(f(x+z)-f(x)-\bigl
\langle \nabla f(x), z\bigr\rangle \mathbf{1}_{\{|z|\le1\}} \bigr)
\frac{C_{d,\alpha
}}{|z|^{d+\alpha}} \,\mathrm{d} z
\\
&&{}+\int_{ \{|z|\le a{|x-y|} \}} \bigl(f(x+z)-f(x)-\bigl\langle \nabla
f(x), z\bigr\rangle \mathbf{1}_{\{|z|\le1\}} \bigr) \biggl(\nu (\mathrm{d}z)-
\frac
{C_{d,\alpha}}{|z|^{d+\alpha}}\, \mathrm{d}z \biggr)
\\
&&{}+\int_{ \{|z|> a{|x-y|} \}} \bigl(f(x+z)-f(x)-\bigl\langle \nabla
f(x), z\bigr\rangle \mathbf{1}_{\{|z|\le1\}} \bigr) \nu(\mathrm {d}z)
\\
&&{}+\bigl\langle b(x),\nabla f(x)\bigr\rangle
\\
&=& Lf(x).
\end{eqnarray*}
This completes the proof.
\end{pf}

\begin{remark}
(1) Here, we give an interpretation of the
construction of the coupling operator $ \widetilde{L} $ above.
If $|x-y|> L_0$, we use the synchronous coupling. If
$|x-y|\le L_0$, then the coupling operator $ \widetilde{L}$
constructed above consists of two parts. Fix any $x,y\in\mathbb
{R}^d$. If
$|z|\le a|x-y|$, then we adopt the coupling by reflection by making
full use of the rotationally invariant measure
$\frac{C_{d,\alpha}}{|z|^{d+\alpha}}\, \mathrm{d}z$; while for the remainder
term, we use the synchronous coupling again, where the
components maintain at each step the same length of jumps (i.e.,
from $(x,y)$ to $(x+z,y+z)$) with the biggest rate $\nu(\mathrm{d}z)$ when
$|z|>a|x-y|$, and with the rate $\nu(\mathrm{d}
z)-\frac{C_{d,\alpha}}{|z|^{d+\alpha}} \,\mathrm{d}z$ when $|z|\le a|x-y|$.
For the coupling by reflection for Brownian motion and diffusion processes,
we refer to \mbox{\cite{Chenbook,LR,Wangbook}}.

(2) Recently, the coupling property of L\'{e}vy processes
has been developed in \cite{SW,BSW,SSW}. The corresponding property
for Ornstein--Uhlenbeck processes with jumps also has been successfully
studied in
\cite{wang1,SW2}. Unlike L\'{e}vy processes and Ornstein--Uhlenbeck
processes with jumps, it is impossible to write
out an explicit expression for transition functions of
the solution to the SDE (\ref{Ito-SDE}) with general drift
term $b(x)$. This observation indicates that all the approaches in
\cite{SW,BSW,SSW,wang1,SW2} are not efficient in the present
setting. This difficulty will be overcome by constructing proper
coupling operators for the Markov generator corresponding to the
solution of the SDE (\ref{Ito-SDE}), as done in \cite{W2014}. However,
different from \cite{W2014} which deals with the corresponding coupling
property by making full use of large jumps part of L\'{e}vy processes,
here to consider the exponential contractivity of the
associated semigroups $(P_t)_{t\ge0}$ with respect to
Wasserstein distances we need a new construction of the coupling
operator. As seen from Propositions \ref{key1} and \ref{key2} below,
the coupling for small jumps part of L\'{e}vy processes [i.e., the
coupling by reflection as mentioned in (1)] is key for our purpose.
\end{remark}

\subsection{Coupling process}\label{section22} In this part, we will
construct a coupling process associated with
the coupling operator $ \widetilde{L}$. For this, we will frequently
talk about the martingale problem for the operator $L$ given by
(\ref{ope-0}) and the coupling operator $ \widetilde{L}$. Let
$\mathcal{D}([0,\infty);\mathbb{R}^d)$ be the space of right continuous
$\mathbb{R}^d$-valued functions having left limits on $[0,\infty)$, equipped
with the Skorokhod topology. For $t\ge0$, denote by $X_t$ the
projection coordinate map on $\mathcal{D}([0,\infty);\mathbb{R}^d)$. A
probability measure $\mathbb{P}^x$ on the Skorokhod space
$\mathcal{D}([0,\infty);\mathbb{R}^d)$ is said to be a solution to the
martingale problem for $(L,C_c^2(\mathbb{R}^d))$ with initial value
$x\in\mathbb{R}^d$ if $\mathbb{P}^x(X_0=x)=1$ and for every $f\in
C_c^2(\mathbb{R}^d)$
\[
\biggl\{f(X_t)-f(x)-\int_0^t
Lf(X_s) \,\mathrm{d}s, t\ge0 \biggr\}
\]
is a $\mathbb{P}^x$-martingale. The martingale problem for
$(L,C_c^2(\mathbb{R}^d))$
is said to be well-posed if it has a unique solution for every
initial value $x\in\mathbb{R}^d$. Similarly, we can define a solution
to the
martingale problem for the coupling operator $ \widetilde{L}$ on
$C_c^2(\mathbb{R}^{2d})$. Note that, in \cite{LM} an equivalence is proved
between the existence of weak solutions to SDEs with jumps and the
existence of solutions to the corresponding martingale problem, by
using a martingale representation theorem. Recently, Kurtz \cite{KU}
studied equivalence between the uniqueness (in sense of
distribution) of weak solutions to a class of SDEs driven by Poisson
random measures and the well-posed solution to martingale problems
for a class of non-local operators using a non-constructive
approach. Note that in our setting the SDE (\ref{Ito-SDE}) has the
pathwise unique
strong solution. According to \cite{BLP}, Theorem~1, page~2, the weak
solution to the SDE (\ref{Ito-SDE}) enjoys the unique (in sense of
distribution) weak solution. This, along with \cite{KU}, Corollary~2.5,
yields that the martingale problem for
$(L,C_c^2(\mathbb{R}^d))$ is well posed.

Let $L_0, a$ be the constants in the definition of the coupling
operator $ \widetilde{L}$. For any $x$, $y\in\mathbb{R}^d$ and $A\in
\mathcal{B}(\mathbb{R}^{2d})$, set
\begin{eqnarray*}
{\mu}(x,y,A)&:=& \frac{1}{2}\int_{ \{(z,\varphi_{x,y}(z))\in
A,|z|\le a{|x-y|}, |x-y|\le L_0 \}} \frac{C_{d,\alpha}}{|z|^{d+\alpha}}
\, \mathrm{d}z
\\
&&{}+\frac{1}{2}\int_{ \{(\varphi_{x,y}(z),z)\in
A,|z|\le a{|x-y|}, |x-y|\le L_0 \}}\frac{C_{d,\alpha
}}{|z|^{d+\alpha}} \,
\mathrm{d}z
\\
&&{}+\int_{ \{(z,z)\in A,|z|\le a{|x-y|}, |x-y|\le
L_0 \}} \biggl(\nu(\mathrm{d}z)-
\frac{C_{d,\alpha}}{|z|^{d+\alpha}} \,\mathrm{d}z \biggr)
\\
&&{}+\int_{ \{(z,z)\in A,|z|>a{|x-y|}, |x-y| \le
L_0
\}\cup \{(z,z)\in A, |x-y|> L_0 \}} \nu(\mathrm{d}z).
\end{eqnarray*}
Then, for any $x$, $y\in\mathbb{R}^d$ and $f\in C^2_b(\mathbb
{R}^{2d})$, we have
\begin{eqnarray*}
&&\widetilde{L} f(x,y)
\\
&&\quad=\int_{\mathbb{R}^{2d}} \bigl[f \bigl((x,y)+(u_1,u_2)
\bigr)-f(x,y)
\\
&&\hspace*{22pt}\qquad{}- \bigl\langle \bigl(\nabla_xf(x,y),\nabla_yf(x,y)
\bigr), (u_1,u_2) \bigr\rangle\mathbf{1}_{\{|u_1|\le1, |u_2|\le1\}}
\bigr] {\mu }(x,y,\mathrm{d}u_1,\mathrm{d}u_2)
\\
&&\qquad {}+\bigl\langle\nabla_x f(x,y),b(x)\bigr\rangle+ \bigl\langle
\nabla_y f(x,y),b(y)\bigr\rangle.
\end{eqnarray*}

Furthermore, for any $h\in C_b(\mathbb{R}^{2d})$, by (A2),
\begin{eqnarray*}
&&\int_{\mathbb{R}^{2d}} h(u)\frac{|u|^2}{1+|u|^2} {\mu}(x,y,\mathrm{d}u)
\\
&&\quad=\int_{ \{|z|\le a{|x-y|}, |x-y|\le L_0 \}}h\bigl(z,\varphi _{x,y}(z)\bigr)
\frac{|z|^2}{1+2|z|^2} \frac{C_{d,\alpha}}{|z|^{d+\alpha}} \,\mathrm{d}z
\\
&&\qquad{}+\int_{ \{|z|\le a{|x-y|}, |x-y|\le L_0 \}
}h\bigl(\varphi_{x,y}(z),z
\bigr)\frac{|z|^2}{1+2|z|^2} \frac{C_{d,\alpha}}{|z|^{d+\alpha}} \,\mathrm{d}z
\\
&&\qquad{}+2\int_{ \{|z|\le a{|x-y|}, |x-y|\le L_0 \}
}h(z,z)\frac{|z|^2}{1+2|z|^2} \biggl(\nu(
\mathrm{d}z)-\frac{C_{d,\alpha}}{|z|^{d+\alpha}} \,\mathrm{d}z \biggr)
\\
&&\qquad {}+2\int_{ \{|z|> a{|x-y|}, |x-y|\le L_0 \} \cup
 \{|x-y|> L_0 \} }h(z,z)\frac{|z|^2}{1+2|z|^2} \nu(
\mathrm{d}z),
\end{eqnarray*}
which implies that $(x,y)\mapsto\int
h(u)\frac{|u|^2}{1+|u|^2} {\mu}(x,y,\mathrm{d}u)$ is a continuous function
on $\mathbb{R}^{2d}$. Note that $b(x)$ is a continuous function on
$\mathbb{R}^d$.
According to \cite{St}, Theorem~2.2, there is a solution to the
martingale problem for $\widetilde{L}$, that is, there exist a probability
space $(\widetilde{\Omega}, \widetilde{\mathcal{F}},
(\widetilde{\mathcal{F}}_t)_{t\ge0}, \widetilde{\mathbb{P}})$ and an
$\bar{\mathbb{R}}^{2d}$-valued process $(\widetilde{X}_t)_{t\ge0}$
such that
$(\widetilde{X}_t)_{t\ge0}$ is
$(\widetilde{\mathcal{F}}_t)_{t\ge0}$-progressively measurable, and
for every $f\in C_b^2(\mathbb{R}^{2d})$,
\[
\biggl\{f(\widetilde{X}_t)-\int_0^{t\wedge e}
\widetilde {L}f(\widetilde {X}_u) \,\mathrm{d}u, t\ge0 \biggr\}
\]
is an $(\widetilde{\mathcal{F}}_t)_{t\ge0}$-local martingale, where
$e$ is the explosion time of $(\widetilde{X}_t)_{t\ge0}$, that is,
\[
e=\lim_{n\to\infty}\inf \bigl\{t\ge0\dvt  |\widetilde{X}_t|\ge n
\bigr\}.
\]
Let\vspace*{1pt} $(\widetilde{X}_t)_{t\ge0}:=(X_t',X_t'')_{t\ge0}$. Then
$(X_t')_{t\ge0}$ and $(X_t'')_{t\ge0}$ are two stochastic processes
on $\mathbb{R}^{d}$. Since $\widetilde{L}$ is the coupling operator
of $L$,
the generator of each marginal process $(X_t')_{t\ge0}$ and
$(X_t'')_{t\ge0}$ is just the operator $L$, and hence both
distributions of the processes $(X_t')_{t\ge0}$ and
$(X_t'')_{t\ge0}$ are solutions to the martingale problem of $L$. In
particular, by our assumption and the remark in the beginning of
this subsection,\vspace*{1.5pt} the processes $(X_t')_{t\ge0}$ and
$(X_t'')_{t\ge0}$ are non-explosive, hence one has $e=\infty$ a.s.
Therefore, the coupling operator $\widetilde{L}$ generates a
non-explosive process $(\widetilde{X}_t)_{t\ge0}$.

Let $T$ be the coupling time of $(X_t')_{t\ge0}$ and
$(X_t'')_{t\ge0}$, that is,
\[
T=\inf\bigl\{t\ge0 \dvt  X_t'=X_t''
\bigr\}.
\]
Then $T$ is an $(\widetilde{\mathcal{F}}_t)_{t\geq0}$-stopping
time. Define a new process $(Y'_t)_{t\ge0}$ as follows:
\[
Y_t'= %
\cases{ X_t'',
& \quad $t< T$;\vspace*{3pt}
\cr
X_t', & \quad $t\ge T$.}
\]
For any $f\in C_b^2(\mathbb{R}^d)$ and $t>0$,
\begin{eqnarray*}
f\bigl(Y_t'\bigr)-\int_0^t
Lf\bigl(Y_s'\bigr) \,\mathrm{d}s&=& f
\bigl(Y_{t\wedge T}'\bigr)-\int_0^{t\wedge T}
Lf\bigl(Y_s'\bigr) \,\mathrm{d}s
\\
&&{} +f\bigl(Y_t'\bigr)-f\bigl(Y_{t\wedge T}'
\bigr)-\int_{t\wedge T}^t Lf\bigl(Y_s'
\bigr) \,\mathrm{d}s
\\
&=& f\bigl(X_{t\wedge T}''\bigr)-\int
_0^{t\wedge T} Lf\bigl(X_s''
\bigr) \,\mathrm{d}s
\\
&&{} +f\bigl(X_t'\bigr)-f\bigl(X_{t\wedge T}'
\bigr)-\int_{t\wedge T}^t Lf\bigl(X_s'
\bigr) \,\mathrm{d}s
\\
&=:&  M_t^1+M_t^2.
\end{eqnarray*}
By the optimal stopping theorem and the facts that both $(X_t')_{t\ge
0}$ and
$(X_t'')_{t\ge0}$ are
solutions to the martingale problem of $L$, $(M_t^1)_{t\ge0}$ and
$(M_t^2)_{t\ge0}$ are martingales and so is $(Y_t')_{t\ge0}$ (see,
e.g., \cite{PW}, Section~3.1, page~251). Since the martingale problem
for the operator $L$ is well-posed, $(Y_t')_{t\ge0}$ and
$(X_t'')_{t\ge0}$ are equal in the distribution. Therefore, we
conclude that
$(X_t',Y_t')_{t\ge0}$ is also a non-explosive coupling process of
$(X_t)_{t\ge0}$ such that $X_t'=Y_t'$ for any $t\ge T$ and
the generator of $(X_t',Y_t')_{t\ge0}$ before the coupling time $T$
is just the coupling operator $\widetilde{L}$. In particular, according to
\cite{CL}, Lemma~2.1, we know
that for any $x$, $y\in\mathbb{R}^d$ and $f\in B_b(\mathbb{R}^d)$,
\[
P_t f(x)={\mathbb{E}^xf\bigl(X_t'
\bigr)}=\widetilde{{\mathbb{E}}}^{(x,y)}f\bigl(X_t'
\bigr)
\]
and
\[
P_t f(y)={\mathbb{E}^yf\bigl(Y_t'
\bigr)}=\widetilde{{\mathbb{E}}}^{(x,y)}f\bigl(Y'_t
\bigr),
\]
where $\widetilde{\mathbb{E}}^{(x,y)}$ is the expectation of the process
$(X_t',Y_t')_{t\ge0}$ with starting point $(x,y)$.

\section{Proofs}\label{section3}
\subsection{Key estimates}
We first assume that $\alpha\in(1,2)$. For any $r>0$, define
\[
\psi(r):= %
\cases{ 1-\mathrm{e}^{-c_1r}, & \quad $r\in[0,2L_0]$;
\vspace*{3pt}\cr
A\mathrm{e}^{c_2(r-2L_0)}+B(r-2L_0)^2+ \bigl(1-\mathrm{e}^{-2c_1L_0}-A
\bigr), & \quad $r\in [2L_0,\infty)$,}
\]
where
\[
A= \frac{c_1}{c_2}\mathrm{e}^{-2L_0c_1},\qquad B= -\frac{(c_1+c_2)c_1}{2}\mathrm{e}^{-2L_0c_1},
\]
$c_2$ is a positive constant such that $c_2\ge20c_1$, that is,
\[
\log\frac{2(c_1+c_2)}{c_2}\le\log2.1,
\]
and $c_1
$ is a positive constant determined by later. With the choice of the
constants $A$ and $B$ above, it is easy to see that $\psi\in
C^2([0,\infty))$. Then we have:

\begin{proposition}\label{key1} Assume that $\alpha\in(1,2)$. Then
there exists a constant $\lambda>0$ such that
for any $x,y\in\mathbb{R}^d$,
\[
\widetilde{L}\psi\bigl(|x-y|\bigr)\le-\lambda\psi\bigl(|x-y|\bigr).
\]
\end{proposition}

\begin{pf} (1) In this part, we treat the case that $x,y\in\mathbb
{R}^d$ with
$|x-y|\le L_0$. First, for any $x, y,z\in\mathbb{R}^d$, by (A3),
\[
\bigl\langle x-y,z+\varphi_{x,y}(z)\bigr\rangle=0
\]
and so
\[
\bigl\langle\nabla_x\psi\bigl(|x-y|\bigr), z+\varphi_{x,y}(z)\bigr
\rangle=0 \quad\mbox{and}\quad \bigl\langle\nabla_y\psi\bigl(|x-y|\bigr), z+
\varphi_{x,y}(z)\bigr\rangle=0.
\]
Therefore,
\begin{eqnarray*}
&& \widetilde{L}\psi\bigl(|x-y|\bigr)
\\
&&\qquad=\frac{1}{2} \biggl[\int_{ \{|z|\le a{|x-y|} \}} \bigl(\psi \bigl( \bigl|x-y+
\bigl(z-\varphi_{x,y}(z)\bigr) \bigr| \bigr)+\psi \bigl( \bigl|x-y-\bigl(z-
\varphi_{x,y}(z)\bigr) \bigr| \bigr)
\\
&&\hspace*{2pt}\quad\qquad\quad{}
-2\psi\bigl(|x-y|\bigr) \bigr)\frac{C_{d,\alpha}}{|z|^{d+\alpha}} \,\mathrm{d} z \biggr]
\\
&&\qquad\quad{} +\psi'\bigl(|x-y|\bigr) \frac{\langle b(x)-b(y),x-y\rangle}{|x-y|}.
\end{eqnarray*}

It is easy to see that $\psi\in C^3([0,2L_0))$ such that $\psi'>0$,
$\psi''<0$ and $\psi'''>0$ on $[0,2L_0)$. Then, for any $0\le
\delta< r \le L_0$,
\begin{eqnarray*}
&&\psi(r+\delta)+\psi(r-\delta)-2\psi(r)=\int_r^{r+\delta}
\,\mathrm{d}s\int_{s-\delta}^s \psi''(u)
\,\mathrm{d}u\leq\psi''(r+\delta)
\delta^2,
\end{eqnarray*}
where in the inequality we have used the fact that $\psi'''>0$ on
$[0,2L_0)$. Hence, according to the definition of $\varphi_{x,y}(z)$
and the inequality above, for all $x,y,z\in\mathbb{R}^d$ with
$|x-y|\le L_0$
and $|z|\le a|x-y|$ with $a\in(0,1/2)$, we have
%
\begin{eqnarray}
&&\psi \bigl( \bigl|x-y+\bigl(z-
\varphi_{x,y}(z)\bigr) \bigr| \bigr)+\psi \bigl( \bigl|x-y-\bigl(z-
\varphi_{x,y}(z)\bigr) \bigr| \bigr)-2\psi\bigl(|x-y|\bigr)
\nonumber
\\
\label{proof-one}
&&\quad=\psi \biggl(|x-y|+\frac{2\langle x-y, z\rangle}{|x-y|} \biggr)+\psi \biggl(|x-y|-
\frac{2\langle x-y, z\rangle}{|x-y|} \biggr)-2\psi\bigl(|x-y|\bigr)
\\
\nonumber
&&\quad\leq4\psi''\bigl((1+2a)|x-y|\bigr)
\frac{\langle x-y,z\rangle^2}{|x-y|^2}.
\end{eqnarray}
Then we deduce that for any $x$, $y\in\mathbb{R}^d$ with $|x-y|\le L_0$,
%
\begin{eqnarray}
\widetilde{L}\psi\bigl(|x-y|\bigr) &\le & 2
\psi''\bigl((1+2a)|x-y|\bigr)\int_{ \{|z|\le a{|x-y|} \}}
\frac{|\langle
x-y,z\rangle|^2}{|x-y|^2}\frac{C_{d,\alpha}}{|z|^{d+\alpha}} \,\mathrm{d}z
\nonumber\\
\nonumber
&&{} + \psi'\bigl(|x-y|\bigr) \frac{\langle
b(x)-b(y),x-y\rangle}{|x-y|}
\\
\nonumber
&=& 2\psi''\bigl((1+2a)|x-y|\bigr)\int
_{ \{|z|\le a{|x-y|} \}} {|z_1|^2}\frac{C_{d,\alpha}}{|z|^{d+\alpha}} \,
\mathrm{d}z
\\
\label{proof-two} %
&&{} + \psi'\bigl(|x-y|\bigr) \frac{\langle
b(x)-b(y),x-y\rangle}{|x-y|}
\\
\nonumber
&=& \frac{2C_{d,\alpha}}{d}\psi''\bigl((1+2a)|x-y|\bigr)
\int_{ \{|z|\le
a{|x-y|} \}} {|z|^2}\frac{1}{|z|^{d+\alpha}} \,\mathrm{d}z
\\
\nonumber
&&{}+ \psi'\bigl(|x-y|\bigr) \frac{\langle
b(x)-b(y),x-y\rangle}{|x-y|}
\\
\nonumber
&\le & \biggl[ - \frac{2C_{d,\alpha} \omega_dL_0^{1-\alpha
}}{d(2-\alpha
)}c_1a^{2-\alpha}
\mathrm{e}^{-2c_1a L_0} +K_1 \biggr] c_1\mathrm{e}^{-c_1|x-y|}|x-y|,
\end{eqnarray}
where in the inequality $z_1$ denotes the first coordinate of $z$, that
is, $z=(z_1,z_2, \ldots, z_d)$, both equalities above follow from the
rotationally invariant property of the measure
$\frac{C_{d,\alpha}}{|z|^{d+\alpha}} \,\mathrm{d}z$, and in the last
inequality we have used (\ref{main-result-1})
and the fact that $\alpha>1$.

Now, taking
\[
C=\frac{2C_{d,\alpha} \omega_dL_0^{1-\alpha}}{d(2-\alpha)},\qquad c_1=(2K_1/C)^{1/(\alpha-1)}
\mathrm{e}^{2L_0/(\alpha-1)}+2, \qquad a=1/c_1,
\]
we find that for any $x$, $y\in\mathbb{R}^d$ with $|x-y|\le L_0$,
\[
\widetilde{L}\psi\bigl(|x-y|\bigr) \le-\frac{C}{2}c_1^{\alpha}
\mathrm{e}^{-2L_0} \mathrm{e}^{-c_1|x-y|}|x-y|.
\]
Since $\psi(0)=0$ and $\psi''\le0$ on $[0,2L_0)$,
\[
\psi(r)\le\psi'(r) r= c_1\mathrm{e}^{-c_1 r} r,\qquad r
\in[0,L_0],
\]
which along with the estimate above yields that
for any $x$, $y\in\mathbb{R}^d$ with $|x-y|\le L_0$,
\[
\widetilde{L}\psi\bigl(|x-y|\bigr) \le-\lambda_1\psi\bigl(|x-y|\bigr),
\]
where $\lambda_1= {C}c_1^{\alpha-1} \mathrm{e}^{-2L_0}/2$.

(2) Second, we consider the case that $x,y\in\mathbb{R}^d$ with $|x-y|>
L_0$. For any $x$, $y\in\mathbb{R}^d$ with $L_0<
|x-y|\le2L_0$, by (\ref{main-result-1}) and the fact that $\psi'>0$,
\begin{eqnarray*}
&& \widetilde{L}\psi\bigl(|x-y|\bigr) \le- c_1 K_2\mathrm{e}^{-c_1|x-y|}|x-y|^{\theta-1}
\le-c_1 K_2 L_0^{\theta-2}
\mathrm{e}^{-c_1|x-y|}|x-y|.
\end{eqnarray*}
On the other hand, also by (\ref{main-result-1}) and the fact that
$\psi
'>0$, for any
$x,y\in\mathbb{R}^d$ with $|x-y|\ge2L_0$,
\begin{eqnarray*}
&& \widetilde{L}\psi\bigl(|x-y|\bigr) \le- K_2 \bigl[Ac_2\mathrm{e}^{c_2(|x-y|-2L_0)}+2B\bigl(|x-y|-2L_0\bigr)
\bigr] |x-y|^{\theta-1}.
\end{eqnarray*}
Next, we consider the function
\[
g(r)=\tfrac{1}{2}Ac_2\mathrm{e}^{c_2(r-2L_0)}+2B(r-2L_0)
\]
on $[2L_0,\infty)$. It
is easy to see that due to the definitions of the constants $A$ and
$B$, there is a unique $r_1\in[2L_0,\infty)$ such
that $g'(r_1)=0$ and
\[
g(r_1)= \frac{-2B}{c_2} \biggl[1-\log\frac{-4B}{Ac_2^2} \biggr]=
\frac
{-2B}{c_2} \biggl[1-\log \frac{2(c_1+c_2)}{c_2} \biggr].
\]
Since
$c_2>0$ is large enough such that
\[
\log \frac{2(c_1+c_2)}{c_2}\le\log2.1,
\]
we have $g(r_1)>0$, which
implies that $g(r)>0$ for all $r\in[2L_0,\infty)$. In
particular,
\[
\tfrac{1}{2}Ac_2\mathrm{e}^{c_2(|x-y|-2L_0)}+2B\bigl(|x-y|-2L_0\bigr)
\ge0
\]
for
any $x,y\in\mathbb{R}^d$ with $|x-y|\ge2L_0$. That is, for any
$x,y\in\mathbb{R}^d$ with
$|x-y|\ge2L_0$,
\begin{eqnarray*}
&& \widetilde{L}\psi\bigl(|x-y|\bigr) \le- \tfrac{1}{2}K_2Ac_2\mathrm{e}^{c_2(|x-y|-2L_0)}|x-y|^{\theta-1}
\le 2^{\theta-3} K_2Ac_2L_0^{\theta-2}\mathrm{e}^{c_2(|x-y|-2L_0)}|x-y|.
\end{eqnarray*}

Combining both estimates above with the definition of $\psi$, we
finally conclude that there is a constant $\lambda_2>0$ such that for
any $x$, $y\in\mathbb{R}^d$ with $|x-y|> L_0$,
\[
\widetilde{L}\psi\bigl(|x-y|\bigr) \le-\lambda_2\psi\bigl(|x-y|\bigr).
\]
This along
with the conclusion of part (1) yields the desired assertion.
\end{pf}

Next, we turn to the case of $\alpha\in(0,1]$. For this, we first take
the constant $a=\frac{1}{4}$ in the definition of the coupling operator
$\widetilde{L}$, and then
change the test function $\psi$ as follows, which is different from
that in the case $\alpha\in(1,2)$. For
any $r>0$, we\vspace*{-2pt} define
\[
\psi(r):= %
\cases{ r-c r^{1+\alpha}, &\quad $r\in[0,2L_0]$;
\vspace*{3pt}\cr
A\mathrm{e}^{c_0(r-2L_0)}+B(r-2L_0)^2+ \bigl(2L_0-c(2L_0)^{1+\alpha}-A
\bigr), &\quad $r\in [2L_0,\infty)$,}
\]
where\vspace*{-3pt}
\[
c=\frac{1}{2^{1+\alpha}(1+\alpha)L_0^{\alpha}},\qquad A= \frac
{1}{2c_0}, \qquad B= -\frac{1}{2} \biggl[
\frac{\alpha}{4L_0}+\frac{c_0}{2} \biggr],\qquad c_0=
\frac{10\alpha}{L_0}.
\]
Due to the choice of the constants above, $\psi\in C^2([0,\infty))$ and
$\psi'(r)>0$ for all $r>0$.

\begin{proposition}\label{key2}
Assume that $\alpha\in(0,1]$.\vspace*{-2pt} If
%
\begin{equation}
\label{con-d1}
\frac{\alpha C_{d,\alpha}\omega
_d3^{\alpha-1}}{8(2-\alpha)d }> K_1 L_0^\alpha,
\end{equation}
then there exists a constant $\lambda>0$ such that
for any $x,y\in\mathbb{R}^d$ with\vspace*{-2pt} $x\neq y$,
\[
\widetilde{L}\psi\bigl(|x-y|\bigr)\le-\lambda\psi\bigl(|x-y|\bigr).
\]
\end{proposition}

\begin{pf}
We mainly follow the proof of Proposition~\ref{key1},
and here we only present the main different steps. For $x,y\in\mathbb{R}^d$
with $|x-y|\le L_0$,\vspace*{-2pt} we have
\begin{eqnarray*}
&& \widetilde{L}\psi\bigl(|x-y|\bigr)
\\[-2pt]
&&\quad=\frac{1}{2} \biggl[\int_{ \{|z|\le ({1}/{4}){|x-y|} \}} \bigl(\psi \bigl(
\bigl|x-y+\bigl(z-\varphi_{x,y}(z)\bigr) \bigr| \bigr)+\psi \bigl( \bigl|x-y-\bigl(z-
\varphi_{x,y}(z)\bigr) \bigr| \bigr)
\\[-2pt]
&&\hspace*{12pt}\qquad{}-2\psi\bigl(|x-y|\bigr) \bigr)\frac{C_{d,\alpha}}{|z|^{d+\alpha}} \,\mathrm{d} z \biggr]
\\[-2pt]
&&\qquad {}+\psi'\bigl(|x-y|\bigr) \frac{\langle b(x)-b(y),x-y\rangle}{|x-y|}.
\end{eqnarray*}
Since $\psi\in C^3((0,2L_0))$ such that $\psi'>0$,
$\psi''<0$ and $\psi'''>0$ on $(0,2L_0)$, one can follow the proof of
(\ref{proof-one}), and get that for all $x,y,z\in\mathbb{R}^d$ with
$0<|x-y|\le
L_0$ and\vspace*{-2pt} $|z|\le
\frac{1}{4}|x-y|$,
\begin{eqnarray*}
&&\psi \bigl( \bigl|x-y+\bigl(z-\varphi_{x,y}(z)\bigr) \bigr| \bigr)+\psi \bigl(
\bigl|x-y-\bigl(z-\varphi_{x,y}(z)\bigr) \bigr| \bigr)-2\psi\bigl(|x-y|\bigr)
\\[-2pt]
&&\quad\leq4\psi'' \biggl(\frac{3}{2}|x-y| \biggr)
\frac{\langle x-y,z\rangle
^2}{|x-y|^2}.
\end{eqnarray*}
Then we follow the argument of (\ref{proof-two}) and deduce that for
any $x$, $y\in\mathbb{R}^d$ with $0<|x-y|\le L_0$,
\begin{eqnarray*}
&& \widetilde{L}\psi\bigl(|x-y|\bigr)\le \biggl[-\frac{\alpha C_{d,\alpha}\omega
_d3^{\alpha-1}}{8(2-\alpha)d L_0^\alpha}+K_1
\biggr]|x-y|.
\end{eqnarray*}
By assumption (\ref{con-d1}), we know that for all $x,y\in\mathbb
{R}^d$ with
$0<|x-y|\le L_0$,
\begin{eqnarray*}
\widetilde{L}\psi\bigl(|x-y|\bigr) &\le & - \biggl(\frac{\alpha C_{d,\alpha}\omega
_d3^{\alpha-1}}{8(2-\alpha)d L_0^\alpha}-K_1
\biggr)|x-y|\le- \biggl(\frac
{\alpha C_{d,\alpha}\omega_d3^{\alpha-1}}{8(2-\alpha)d L_0^\alpha
}-K_1 \biggr)\psi\bigl(|x-y|\bigr)\\
&=: & -
\lambda_1\psi\bigl(|x-y|\bigr).
\end{eqnarray*}

Next, we turn to the case that $x,y\in\mathbb{R}^d$ with $|x-y|>
L_0$. For any
$x$, $y\in\mathbb{R}^d$ with $L_0< |x-y|\le
2L_0$, by (\ref{main-result-1}) and $\psi'>0$,
\begin{eqnarray*}
\widetilde{L}\psi\bigl(|x-y|\bigr) &\le & - K_2 \bigl(1-c(1+\alpha)|x-y|^\alpha
\bigr)|x-y|^{\theta-1}
\\
&\le & -K_2 L_0^{\theta-2} \biggl(1-
\frac{1}{2^{1+\alpha}L_0^\alpha
}|x-y|^\alpha \biggr)|x-y|.
\end{eqnarray*}
On the other hand, also by (\ref{main-result-1}) and the fact that
$\psi
'>0$, for any
$x,y\in\mathbb{R}^d$ with $|x-y|\ge2L_0$,
\begin{eqnarray*}
\widetilde{L}\psi\bigl(|x-y|\bigr) &\le & - K_2 \bigl[Ac_0\mathrm{e}^{c_0(|x-y|-2L_0)}+2B\bigl(|x-y|-2L_0\bigr)
\bigr] |x-y|^{\theta-1}.
\end{eqnarray*}
Now, we consider again the function
\[
g(r)=\tfrac{1}{2}Ac_0\mathrm{e}^{c_0(r-2L_0)}+2B(r-2L_0)
\]
on $[2L_0,\infty)$. It
is easy to see that due to the definitions of the constants $A$ and
$B$, there is a unique $r_1\in[2L_0,\infty)$ such
that $g'(r_1)=0$ and
\[
g(r_1)= \frac{-2B}{c_0} \biggl[1-\log\frac{-4B}{Ac_0^2} \biggr]=
\frac
{-2B}{c_0} \biggl[1-\log \biggl(2+\frac{\alpha}{L_0c_0} \biggr) \biggr].
\]
Noticing that $c_0={10\alpha}{L_0}^{-1}$, we get
\[
\log \biggl(2+\frac{\alpha}{L_0c_0} \biggr)= \log2.1,
\]
and so $g(r_1)>0$, which
implies that $g(r)>0$ for all $r\in[2L_0,\infty)$. In
particular,
\[
\tfrac{1}{2}Ac_0\mathrm{e}^{c_0(|x-y|-2L_0)}+2B\bigl(|x-y|-2L_0\bigr)
\ge0
\]
for
any $x,y\in\mathbb{R}^d$ with $|x-y|\ge2L_0$. That is, for any
$x,y\in\mathbb{R}^d$ with
$|x-y|\ge2L_0$,
\begin{eqnarray*}
\widetilde{L}\psi\bigl(|x-y|\bigr) &\le & - \tfrac{1}{2}K_2Ac_0\mathrm{e}^{c_0(|x-y|-2L_0)}|x-y|^{\theta-1}.
\end{eqnarray*}

According to both estimates above and the definition of $\psi$, we
finally conclude that there is a constant $\lambda_2>0$ such that for
any $x$, $y\in\mathbb{R}^d$ with $|x-y|> L_0$,
\[
\widetilde{L}\psi\bigl(|x-y|\bigr) \le-\lambda_2\psi\bigl(|x-y|\bigr).
\]
This along
with the conclusion above yields the desired assertion.
\end{pf}

\begin{remark}\label{main-1-1} According to the argument above, we
can easily improve (\ref{con-d1}), for example, by taking $\psi
(r)=c_1r-c_2r^{1+\alpha'}$ for $r\in[0,2L_0]$ and changing
the integral domain $\{z \dvt  |z|\le\frac{1}{4}|x-y|\}$ in the
definition of the coupling operator $\widetilde{L}$ into $\{z \dvt |z|\le
a|x-y|\}$ with some proper choices of $c_1,c_2>0$, $\alpha'\in
(0,\alpha
]$ and $a\in
(0,1/2)$. For simplicity, here we just set $c_1=1$, $c_2=c$, $\alpha
'=\alpha$ and
$a=1/4$.
\end{remark}

\subsection{Proofs of Theorem~\texorpdfstring{\protect\ref{main-result}}{1.2} and
Example~\texorpdfstring{\protect\ref{1-example-2}}{1.3}}\label{section2}
We divide the proof of Theorem~\ref
{main-result} into two parts.

\begin{pf*}{Proof of Theorem \protect\ref{main-result} for
$|x-y|\le L_0$ or $\theta=2$}
We will make full use of the coupling process $(X'_t,Y'_t)_{t\ge0}$ constructed
in Section~\ref{section22}. Denote by $\widetilde{\mathbb
{P}}^{(x,y)}$ and
$\widetilde{\mathbb{E}}^{(x,y)}$ the distribution and the expectation of
$(X'_t,Y'_t)_{t\ge0}$ starting from $(x,y)$, respectively. For any
$t>0$ set $r_t=|X_t'-Y_t'|$, and for $n\geq1$ define the stopping time
\[
T_n=\inf\bigl\{t>0 \dvt   r_t\notin[1/n, n]\bigr\}.
\]

For any $x$, $y\in\mathbb{R}^d$ with $|x-y|>0$,
we take $n$ large enough such that $1/n<|x-y|<n$. Let $\psi$ be the
function given in Proposition~\ref{key1} if $\alpha\in(1,2)$ or the
function given in Proposition~\ref{key2} if $\alpha\in(0,1]$. Then
\begin{eqnarray*}
&& \widetilde{\mathbb{E}}^{(x,y)}\psi \bigl(\bigl|X'_{t\wedge
T_{n}}-Y'_{t\wedge
T_{n}}\bigr|
\bigr)
\\
&&\quad=\psi\bigl(|x-y|\bigr)+\widetilde{\mathbb{E}}^{(x,y)} \biggl(\int
_0^{t\wedge T_{n}} \widetilde{L} \psi \bigl(\bigl|X'_{s}-Y'_{s}\bigr|
\bigr) \,\mathrm{d}s \biggr)
\\
&&\quad\le\psi\bigl(|x-y|\bigr)-\lambda\widetilde{\mathbb{E}}^{(x,y)} \biggl(\int
_0^t \psi \bigl(\bigl|X'_{{s\wedge T_{n}}}-Y'_{{s\wedge T_{n}}}\bigr|
\bigr) \,\mathrm{d}s \biggr).
\end{eqnarray*}
Therefore,
\[
\mathbb{E} \bigl[\psi(r_{t\wedge T_n}) \bigr]\leq\psi(r_0)
\mathrm{e}^{-\lambda t}.
\]
Since the coupling process $(X'_t,Y'_t)_{t\ge0}$ is non-explosive, we
have $T_n\uparrow T$ a.s. as $n\to\infty$,
where $T$ is the coupling time of the process $(X_t',Y_t')$. Thus, by
Fatou's lemma, letting $n\to\infty$ in the above inequality
gives us
\[
\mathbb{E} \bigl[\psi(r_{t\wedge T}) \bigr]\leq\psi(r_0)\mathrm{e}^{-\lambda t}.
\]
Thanks to our convention that $Y'_t=X'_t$ for $t\geq T$, we have
$r_t=0$ for all $t\geq T$,
and so
\[
\mathbb{E}\psi(r_t)\leq\psi(r_0)\mathrm{e}^{-\lambda t}.
\]
That is,
\[
\mathbb{E}\psi\bigl(|X_t-Y_t|\bigr)\leq\psi\bigl(|x-y|\bigr)\mathrm{e}^{-\lambda t}.
\]
As a result, if
$|x-y|\le L_0$, then for any $p\ge1$ and $t>0$,
%
\begin{equation}
\label{small} %
\mathbb{E}|X_t-Y_t|^p
\le C(p)\mathbb{E}\psi\bigl(|X_t-Y_t|\bigr) \le C_1
\mathrm{e}^{-\lambda t}|x-y|,
\end{equation}
where the first inequality follows from the definitions of the test
function $\psi$ in Propositions \ref{key1} and \ref{key2}.

Now for any $x,y\in\mathbb{R}^d$ with $|x-y|>L_0$, take
$n:= [|x-y|/L_0 ]+ 1\geq2$. We have
%
\begin{equation}
\label{proof.2} %
\frac{n}{2}\le n-1\leq
\frac{|x-y|}{L_0}\le n.
\end{equation}
Set $x_i=x+i(y-x)/n$ for $i=0,1,\ldots, n$. Then $x_0=x$ and
$x_n=y$; moreover, (\ref{proof.2}) implies
$|x_{i-1}-x_{i}|=|x-y|/n\le L_0$ for all $i=1,2,\ldots,n$.
Therefore, by (\ref{small}) and (\ref{proof.2}),
\begin{eqnarray*}
W_{p}(\delta_x P_t,\delta_y
P_t)&\leq & \sum_{i=1}^n
W_{p}(\delta _{x_{i-1}} P_t,\delta_{x_{i}}
P_t)
\\
&\leq &  C_1^{1/p} \mathrm{e}^{-{\lambda}t/p}\sum
_{i=1}^n |x_{i-1}-x_{i}|^{1/p}
\\
&\leq  & C_1^{1/p} \mathrm{e}^{-{\lambda}t/p} nL_0^{1/p}
\\
&\leq & 2C_1^{1/p}L_0^{1/p-1}\mathrm{e}^{-{\lambda}t/p}
|x-y|
\\
&=: & C_2\mathrm{e}^{-{\lambda t}/p}|x-y|.
\end{eqnarray*}
In particular, the proof of the first assertion for $\theta=2$ in
Theorem~\ref{main-result} is completed. On the other hand, from
(\ref{small}) and the conclusion above, we also get the second
assertion for $\theta>2$ with $|x-y|\le1$ and all $t>0$, or with
$|x-y|>1$ and $0<t\le1$.
\end{pf*}

Next, we turn to:
\begin{pf*}{Proof of Theorem \protect\ref{main-result} for $|x-y|>
L_0$ and $\theta>2$}
For $|x-y|> L_0$, we use the synchronous coupling and
the assertion of Theorem $\ref{main-result}$ for $|x-y|\le L_0$. In
detail, with (\ref{Ito-SDE}), let $(X_t,Y^{(2)}_t)_{t\ge0}$ be the
coupling process on $\mathbb{R}^{2d}$ such that its distribution is
the same as
that of $(X'_t,Y'_t)_{t\ge0}$ constructed
in Section~\ref{section22}. We now consider
%
\begin{equation}
\label{coupling-1}
\mathrm{d}Y_t=
\cases{
\mathrm{d}Z_t+b(Y_t)\, \mathrm{d}t, & \quad $0\leq t<
T_{L_0}$,\vspace*{3pt}
\cr
\mathrm{d}Y_t^{(2)}, &\quad $T_{L_0}\le t< T$,}
\end{equation}
where
\[
T_{L_0}=\inf\bigl\{t>0 \dvt |X_t-Y_t|\le
L_0\bigr\}
\]
and $T=\inf\{t>0  \dvt X_t=Y_t\}$ is the coupling time. For $t\geq T$, we
still set $Y_t=X_t$. Therefore, the difference process
$(D_t)_{t\ge0}:=(X_t-Y_t)_{t\ge0}$ satisfies
%
\[
\label{difference-1} \mathrm{d}D_t=\bigl(b(X_t)-b(Y_t)
\bigr) \,\mathrm{d}t, \qquad t<T_{L_0}.
\]
%
Note that the equality above implies that $t\mapsto D_t$ is a
continuous function on $[0,T_{L_0})$ such that
$\lim_{t\to T_{L_0}-} |D_t|=L_0$.
As a result,
\[
\mathrm{d}|D_t|^2=2 \bigl\langle D_t,b(X_t)-b(Y_t)
\bigr\rangle \,\mathrm{d}t,\qquad t<T_{L_0}.
\]
Still denoting by $r_t=|D_t|$, we get from (\ref{main-result-1})
that
\[
\mathrm{d}r_t\leq-K_2r_t^{\theta-1} \,
\mathrm{d}t, \qquad t<T_{L_0},
\]
which implies that
%
\begin{equation}
\label{rrr2}
T_{L_0}\le\frac{1}{K_2(2-\theta)} \bigl(|x-y|^{2-\theta
}-L_0^{2-\theta
}
\bigr) \le\frac{L_0^{2-\theta}}{K_2(\theta-2)}=: t_0
\end{equation}
since $\theta>2$ and the continuity of $t\mapsto r_t$ on $[0,T_{L_0})$.

Therefore, for any $x,y\in\mathbb{R}^d$ with $|x-y|> L_0$,
$p\ge1$ and $t>t_0$, we have
\begin{eqnarray*}
\mathbb{E}|X_t-Y_t|^p&
=& \mathbb{E} \bigl[\mathbb{E}^{(X_{T_{L_0}},
Y_{T_{L_0}})}|X_{t-T_{L_0}}-Y_{t-T_{L_0}}|^p
\bigr]
\\
& \le  & C_1\mathbb{E} \bigl[|X_{T_{L_0}}-Y_{T_{L_0}}|\mathrm{e}^{-\lambda
(t-T_{L_0})}
\bigr]
\\
&\le &  C_1 L_0\exp(\lambda t_0)
\mathrm{e}^{-\lambda t},
\end{eqnarray*}
where in the first inequality we have used (\ref{small}), and the
last inequality follows from (\ref{rrr2}) and the fact that
$|X_{T_{L_0}}-Y_{T_{L_0}}|\le L_0$. In particular, we have for all
$|x-y|>L_0$ and $t>t_0$,
\[
\mathbb{E}|X_t-Y_t|^p\le C_3
\mathrm{e}^{-\lambda t}.
\]
Combining with all conclusions above, we complete the proof of the
second assertion in Theorem~\ref{main-result}.
\end{pf*}

We finally present the
following.
\begin{pf*}{Proof of Example \protect\ref{1-example-2}}
In this
example,
\[
b(x)=\nabla V(x)=2\beta|x|^{2\beta-2}x.
\]
It follows from the proof of \cite{CG}, Example~5.3,  that for any
$x,y\in
\mathbb{R}^d$,
%
\begin{equation}
\label{cg-1}
\bigl\langle b(x)-b(y),x-y\bigr\rangle\le-\beta
2^{4-3\beta}|x-y|^{2\beta}.
\end{equation}
Then, (\ref{main-result-1}) holds with $K_2=\beta2^{4-3\beta}$,
$\theta=2\beta$ and any positive constants $K_1,L_0$. In particular,
(\ref{main-result-1-1}) holds for all $\alpha\in(0,1]$ and $K_1,
L_0>0$ small enough. Then the required assertion is a direct
consequence of Theorem~\ref{main-result}.
\end{pf*}

\section*{Acknowledgements}
The author would like to thank
Professor Feng-Yu Wang and the referee for helpful comments and careful
corrections.
Financial support through NSFC (No. 11201073), JSPS (No. 26$\cdot$04021),
NSF-Fujian (No. 2015J01003) and the Program for Nonlinear Analysis and Its
Applications (No. IRTL1206)  are gratefully acknowledged.



%



\printhistory
\end{document}